\documentclass[a4paper]{amsart}

\title[Borel rigidity for higher rank lattices]{Borel rigidity for higher rank lattices}

\author[F.~Calderoni]{Filippo Calderoni}
\address{Department of Mathematics, Rutgers University, 
Hill Center for the Mathematical Sciences,
110 Frelinghuysen Rd.,
Piscataway, NJ 08854-8019}
\email{filippo.calderoni@rutgers.edu}

\author[D.~Rupp]{Danae Rupp}
\address{Department of Mathematics, Rutgers University, 
Hill Center for the Mathematical Sciences,
110 Frelinghuysen Rd.,
Piscataway, NJ 08854-8019}
\email{danae.rupp@rutgers.edu}

\date{\today} 

\subjclass[2020]{Primary: 03E15, 37A20.}

\thanks{The first author was supported by the NSF Grant DMS -- 2348819.
We thank Konrad Wr\'obel for pointing out that the proof of Theorem~\ref{thm : containment shift} extends to Gaussian actions. (See Remark~\ref{remark : Gaussian}).}

\keywords{}

\usepackage[T1]{fontenc}
\usepackage[osf]{mathpazo}

\usepackage{eucal}

\usepackage{microtype}

\usepackage[colorlinks=true,linkcolor=cyan,citecolor=cyan]{hyperref}

\usepackage{amssymb,
enumitem,
mathrsfs,
mathtools,
tikz,
upgreek,
verbatim,
hyphenat,
url,
amsmath
}
\usepackage[T1]{fontenc}
\usepackage[shadow]{todonotes}

\newenvironment{enumerate-(a)}{\begin{enumerate}[label={\upshape (\alph*)}, leftmargin=2pc]}{\end{enumerate}}
\newenvironment{enumerate-(a)-r}{\begin{enumerate}[label={\upshape (\alph*)}, leftmargin=2pc,resume]}{\end{enumerate}}
\newenvironment{enumerate-(a)-5}{\begin{enumerate}[label={\upshape (\alph*)}, leftmargin=2pc,start=5]}{\end{enumerate}}
\newenvironment{enumerate-(A)}{\begin{enumerate}[label={\upshape (\Alph*)}, leftmargin=2pc]}{\end{enumerate}}
\newenvironment{enumerate-(A)-r}{\begin{enumerate}[label={\upshape (\Alph*)}, leftmargin=2pc,resume]}{\end{enumerate}}
\newenvironment{enumerate-(i)}{\begin{enumerate}[label={\upshape (\roman*)}, leftmargin=2pc]}{\end{enumerate}}
\newenvironment{enumerate-(i)-r}{\begin{enumerate}[label={\upshape (\roman*)}, leftmargin=2pc,resume]}{\end{enumerate}}
\newenvironment{enumerate-(I)}{\begin{enumerate}[label={\upshape (\Roman*)}, leftmargin=2pc]}{\end{enumerate}}
\newenvironment{enumerate-(I)-r}{\begin{enumerate}[label={\upshape (\Roman*)}, leftmargin=2pc,resume]}{\end{enumerate}}
\newenvironment{enumerate-(1)}{\begin{enumerate}[label={\upshape (\arabic*)}, leftmargin=2pc]}{\end{enumerate}}
\newenvironment{enumerate-(1)-r}{\begin{enumerate}[label={\upshape (\arabic*)}, leftmargin=2pc,resume]}{\end{enumerate}}
\newenvironment{itemizenew}{\begin{itemize}[leftmargin=2pc]}{\end{itemize}}

\newtheorem{theorem}{Theorem}[section]
\newtheorem{lemma}[theorem]{Lemma}
\newtheorem{corollary}[theorem]{Corollary}
\newtheorem{proposition}[theorem]{Proposition}

\newtheorem{claim}{Claim}[theorem]
\theoremstyle{definition}
\newtheorem{definition}[theorem]{Definition}

\newtheorem{example}[theorem]{Example}

\theoremstyle{remark}
\newtheorem{remark}[theorem]{Remark}

\numberwithin{equation}{section}

\DeclareMathOperator{\rank}{rank}

\DeclareMathOperator{\GL}{GL}
\DeclareMathOperator{\SL}{SL}

\DeclareMathOperator{\SO}{SO}

\DeclareMathOperator{\Aut}{Aut}

\DeclareMathOperator{\covol}{covol}
\DeclareMathOperator{\Inn}{Inn}

\begin{document}

\begin{abstract}
In this paper we focus on rigidity in descriptive set theory. Building on cocycle superrigidity results of Popa, Furman and Ioana, we prove rigidity in the purely Borel context for the orbit equivalence relations arising from the shift actions and translation actions of lattices in semisimple Lie groups with higher rank. Our results sharpen work of Adams and Kechris.
\end{abstract}

\maketitle

\section{Introduction}

\subsection{Background}
Rigidity is the phenomenon
when every single element of a class of mathematical objects is determined by a smaller
amount of data that one would expect. The typical scheme of a rigidity theorem shows
that whenever two objects are equivalent up to some notion of isomorphism, then they
are equivalent in a stronger sense.

A large supply of rigidity results has been proven in the context of lattices in higher rank Lie groups. For example, a celebrated result of Mostow and Margulis states that two irreducible lattices in centerless semisimple Lie groups are isomorphic if and only if there is a continuous isomorphism between the ambient Lie groups. (See also Theorem~\ref{thm:mostowmargulis} below.)

This paper focuses on rigidity for group actions and is inspired by the breakthrough work of Zimmer~\cite{Zim}. We recall some standard definitions below.

Suppose that \(\Gamma\) is a discrete countable group and \((X, \mu)\) is a
standard Borel space with probability measure \(\mu\) on which \(\Gamma\) acts in a Borel fashion. We say that the action \(\Gamma\times X \to X, (g,x)\mapsto g\cdot x\) is \emph{probability measure preserving} (\emph{p.m.p.} for short) if
\(\mu(g\cdot A) = \mu(A)\) for all Borel \(A \subseteq X\) and \(g \in \Gamma\). Moreover we are particularly interested in the case when the action of \(\Gamma\) on \(X\) is ergodic, i.e. for all Borel \(A\subseteq X\) and all \(g \in\Gamma\), we have \(\mu(gA\triangle A) = 0\), which implies that \(A\) is null or conull.

Two p.m.p. actions \(\Gamma \curvearrowright (X, \mu)\) and \(\Lambda \curvearrowright (Y, \nu)\) of countable groups \(\Gamma\) and \(\Lambda\)
are \emph{orbit equivalent} if there is a
measure preserving isomorphism \(\theta \colon X \to Y\) such that \(\theta(\Gamma \cdot x) = \Lambda \cdot \theta(x)\), for almost every \(x \in X\).

In one of his seminal papers~\cite{Zim80}, Zimmer proved the following rigidity result.

\begin{theorem}[Zimmer~1980]
\label{thm : Zimmer rigidity}
Suppose that \(G_0\) and \(G_1\) are connected
centerless simple Lie groups of \(\mathbb{R}\)-rank at least \(2\) and that \(\Gamma\) and \(\Delta\) are lattices in \(G_0, G_1\),
respectively. Let \(\Gamma\) and \(\Delta\) act freely on the ergodic probability space \((X_0,\mu_0)\) and \((X_1,\mu_1)\) respectively.

If the actions \(\Gamma \curvearrowright (X_0, \mu_0)\) and \(\Lambda \curvearrowright (X_1, \mu_1)\) are (weakly) orbit equivalent, then \(G_0 \cong G_1\).
\end{theorem}

Consequently, for \(m,n\geq 2\). The group action \(\SL_m(\mathbb{Z}) \curvearrowright \mathbb{T}^m  \) is orbit equivalent to \(\SL_n(\mathbb{Z}) \curvearrowright \mathbb{T}^n  \) only if \(m=n\). This contrasts with a famous anti-rigidity result by Dye~\cite{Dye63} and Ornstein--Weiss~\cite{OrnWei} stating that any two free ergodic p.m.p. actions of two amenable groups are orbit
equivalent.
In particular, the irrational toral rotation actions \(\mathbb{Z}^m \curvearrowright(\mathbb{T}^m,\mu)\) and \(\mathbb{Z}^n \curvearrowright(\mathbb{T}^n,\mu)\) are orbit equivalent for any \(m,n\in\mathbb{N}\).

Moreover, Ioana~\cite{Ioa09,Ioa11Inv} and Epstein~\cite{Eps08} showed that for non-amenable groups, the situation is completely different. 

\begin{theorem}[Ioana~2007; Epstein~2008]
    For every non-amenable group \(\Gamma\), there are uncountably many orbit inequivalent free ergodic \(\Gamma\)-actions.
\end{theorem}

Back to rigidity for higher rank
lattices, the seminal work of Zimmer was fundamental to the discovery by Adams and Kechris that there are continuum many incomparable countable Borel equivalence relations up to Borel equivalence.

 If \(E\) and \(F\) are Borel equivalence relations on the standard Borel spaces \(X\) and \(Y\), respectively, then a Borel map \(f \colon X \to Y\) is said to be a \emph{homomorphism} from \(E\) to \(F\) if for all \(x, y \in X\),
\[
x\mathbin{E}y\implies f(x)\mathbin{F}f(y).
\]
A \emph{weak Borel reduction} from \(E\) to \(F\) is a countable-to-one Borel homomorphism.
Whenever \(f\) satisfies the stronger property that for all \(x, y \in X\),
\[
x\mathbin{E}y \iff f(x)\mathbin{F}f(y),
\]
then
we say that
\(E\) is \emph{Borel reducible} to \(F\), that
\(f\) is a \emph{Borel reduction}, and we write \(E\leq_{B}F\).

We say that a Borel equivalence relation is \emph{countable} if it has countable classes. 
A classical result of Feldman and Moore~\cite{FelMoo}
establishes that all countable Borel equivalence relations arise as the orbit equivalence relation induced by the Borel action of some countable discrete group. Whenever \(\Gamma\) be a countable discrete group acting on the standard Borel space \(X\) in a Borel way, we denote by \(E^X_\Gamma\) the corresponding induced orbit equivalence relation on \(X\), whose equivalence classes are the \(\Gamma\)-orbits. That is, we define \(E^X_\Gamma\) by letting \(x\mathbin{E}^X_Gy \iff \exists g\in G (g\cdot x = y)\), for all \(x,y\in X\).

For example, whenever \(\Gamma\) is a countable group, we may consider the Polish space \(2^\Gamma=\{f\colon \Gamma\to \{0,1\}\}\) with the product topology. We let \(\Gamma\) act on \(2^\Gamma\) by left-shift, so \(g\cdot f(x) = f(g^{-1} x)\), for all \(g\in \Gamma\) and \(f\in 2^\Gamma\). Define the \emph{free part} of the action of the shift-action of \(G\) by
$$(2)^\Gamma=\{x \in 2^\Gamma \mid g\cdot x\ne x \text{ for all }1\ne g \in \Gamma\}.$$
It is easily seen that \((2)^\Gamma\) is a \(\Gamma\)-invariant \(G_\delta\) subset of \(2^\Gamma\). Therefore, we can regard \((2)^\Gamma\) as a standard Borel space and we  denote the equivalence relation induced by the action of \(\Gamma\) on \((2)^\Gamma\) as $F_\Gamma$.

As discussed by Thomas in his ICM address~\cite{Tho06}, the methods developed by Adams and Kechris yield the following Borel rigidity result for higher rank lattices in the purely Borel context: 

\begin{theorem}[Adams--Kechris~2000]
\label{thm : AdaKec rigidity}
    Suppose that \(G_0\) and \(G_1\) are connected
centerless simple Lie groups of \(\mathbb{R}\)-rank at least \(2\) and that \(\Gamma_0\) and \(\Gamma_1\) are lattices in \(G_0, G_1\),
respectively. Let \(\Gamma_0\) and \(\Gamma_1\) act freely on the ergodic probability space \((X_0,\mu_0)\) and \((X_1,\mu_1)\) respectively.
If \(E^{X_0}_{\Gamma_0}\sim_B E^{X_1}_{\Gamma_1}\), then \(G_0 \cong G_1\).
\end{theorem}

\subsection{Content of the paper} In this paper we explore the interplay between Borel equivalence and superrigidity more in depth for specific actions. First, we prove a strengthening of Theorem~\ref{thm : AdaKec rigidity} in the more specific context of free shift actions.

\begin{theorem}
\label{thm:rigidity shift}
    Let $G$ and \(H\) be centerless semisimple Lie groups with \(\mathbb{R}\)-rank at least \(2\), and let $\Gamma$ and \(\Delta\) be irreducible lattices in $G$ and \(H\), respectively. If $F_\Gamma \sim_BF_\Delta$, then $\Gamma\cong\Delta$.
\end{theorem}

Consequently, one gets \(G\cong H\) by Mostow--Margulis rigidity.
We also note that the same conclusion is obtained whenever \(F_\Gamma\) and \(F_\Delta\) are orbit equivalent.

While Theorem~\ref{thm : Zimmer rigidity} relies on Zimmer's cocycle superrigidity, our proof of Theorem~\ref{thm:rigidity shift} uses the more specialized cocycle superrigidity theory for Bernoulli actions due to Popa~\cite{Pop07,Pop08}, which allows us to obtain a stronger conclusion.

As an immediate application, we discuss principal congruence subgroups.
For a prime number \(p\), let \(\Gamma(3,p)\) be the kernel of the surjective map
\(\SL_3(\mathbb{Z}) \to \SL_3( \mathbb{F}_p)\)
induced by reduction modulo \(p\). (Here, \(\mathbb{F}_p\) is the field of \(p\) elements.)
Since the image of the reduction map is finite, \(\Gamma(3,p)\) has finite index in \(\SL_3(\mathbb{Z})\). Therefore, it is an irreducible lattice in \(\SL_3(\mathbb{R})\).

For different prime numbers \(p,q\), it is well known that \(\Gamma(3, p)\not\cong\Gamma(3, q)\). In fact, their abelianizations are not isomorphic. (E.g. see~\cite[Theorem~1.1]{ImoKob} as well as~\cite{LeeSzc})

Therefore,  we obtain the following, which is an easy corollary of Theorem~\ref{thm:rigidity shift}. 

\begin{corollary}
\label{cor : congruence subgroup}
    For any prime numbers \(p, q\), if \(F_{\Gamma(3,p)} \sim_B F_{\Gamma(3,q)}\) then \(p=q\).
\end{corollary}

As both \(\Gamma(3,p)\) and \(\Gamma(3,q)\) are both lattices in \(\SL_3(\mathbb{R})\), Corollary~\ref{cor : congruence subgroup} does not follow from Adams and Kechris's rigidity result.

\begin{remark}
   The argument of Theorem~\ref{thm:rigidity shift} can be run for every orbit equivalence relation induced by a free ergodic p.m.p. malleable actions on non-atomic probability spaces.
\end{remark}

Next, we prove an analog of Theorem~\ref{thm:rigidity shift} for the translation action. Here the acting groups are \(S\)-arithmetic groups, and therefore lattices in a product of locally compact groups, acting by translation on the ``ambient'' compact Lie group.

We briefly recall our notation. 
Let \(k\) be a global field with a unique archimedean place, called \(\infty\). Let \(S\) and \(T\) be a finite set of
places of \(k\), containing \(\infty\) and at least one non-archimedean place.
Whenever \(G\) is an algebraic semisimple group over \(k\) we define
\(G_S = \prod_{s\in S}G(k_s)\).

\begin{theorem}
\label{thm:rigidity translation}
Let \(S,T\) be finite non-empty sets and let \(G,H\) be centerless connected almost \(k\)-simple groups such that:
\begin{enumerate-(1)}
    \item
    \(G(k_\infty)\) and \(H(k_\infty)\) are compact.
    \item For each \(s\in S\) and \(t\in T\), both
    \(\rank_{k_s} G\) and \(\rank_{k_t} H\) are at least \(2\).
    \item 
    \(G(k_s)\) and \(G(k_t)\) have no non-trivial compact factors.
\end{enumerate-(1)}
    Let \(\Gamma\) and \(\Delta\) be irreducible lattices in \(G_S\) and \(H_T\), respectively, and let \(\Gamma\curvearrowright G(k_\infty)\) and \(\Lambda\curvearrowright H(k_\infty)\) be the translation actions. If \(E^{G(k_\infty)}_{\Gamma} \sim_B E^{H(k_\infty)}_{\Delta}\), then \(\Gamma \cong \Delta\).
\end{theorem}

Our results contribute to the literature of rigidity in descriptive set theory~\cite{AdaKec,Cal23,Cos10, Tho03APAL}. See also \cite[Section~6]{Kec24} for a full account of Borel (or set-theoretic) rigidity in descriptive set theory.

In Section~2 we discuss preliminaries. In Section~3 we prove the two main results  of this paper (Theorem~\ref{thm:rigidity shift} and Theorem~\ref{thm:rigidity translation}.) In Section~4 we discuss how this is related to Adams's and Thomas's solutions to the containment problem. In particular, we see how the former follows from Theorem~\ref{thm:rigidity translation}. Since our results might be interesting to a broad audience, we wrote this paper keeping in mind that the readers might very well have different backgrounds. As a result, the preliminary section is substantial.


\section{Preliminaries}

\label{sec : preliminaries}

\subsection{Ergodic theory}
If \(G\) is a locally compact second countabel (lcsc) group and acts on \(X\) in a Borel way, we say that \(X\) is standard Borel \(G\)-space.
Let \(G\) be a lcsc group and \(X\) be a standard Borel \(G\)-space. A standard Borel probability measure \(\mu\) on \(X\) is \emph{\(G\)-invariant} if \(\mu(g A) = \mu(A)\), for every \(g\in G\) and Borel \(A\subseteq X\). We recall the definition of ergodic action. Whenever \(\mu\) is \(G\)-invariant, then the action of \(G\) on \((X,\mu)\) is  \emph{ergodic} if, for every \(G\)-invariant Borel subset \(A\subseteq X\), either \(\mu(A)=0\) or \(\mu(A)=1\). Recall a function \(f\colon X \to Y\) between standard Borel \(G\)-spaces is \emph{\(G\)-invariant} if \(f(g\cdot x) = f(x)\) for all \(g\in G\) and \(x\in X\). The following proposition is well-known. (E.g., see \cite[Theorem~1.3]{BekMay}.)

\begin{proposition}
\label{prop:ergodic}
If \(\mu\) is a \(G\)-invariant probability measure on the standard Borel \(G\)-space \(X\), then the following are equivalent.
\begin{enumerate-(i)}
\item
\label{item:ergodicity(i)}
The action of \(G\) on \((X,\mu)\) is ergodic.
\item
\label{item:ergodicity(ii)}
If \(Y\) is a standard Borel space and \(f\colon X \to Y\) is a \(G\)-invariant Borel function, then \(f\) is essentially constant; i.e., there is some Borel \(Z\subseteq X\) with \(\mu(Z)=1\) such that \(f\restriction Z\) is constant.
\end{enumerate-(i)}
\end{proposition}

If there is a unique $G$-invariant probability measure $\mu$ on $X$, then we say that the action of $G$ on $X$ is \textit{uniquely ergodic}. In this case, it is well-known that \(\mu\) is an ergodic measure. (E.g., see~\cite[Definition~I.3.6]{BekMay}.)

Suppose that \(X\) and \(Y\) are standard Borel spaces and that \(f\colon X\to Y\) is a Borel function. Then, for every probability measure \(\nu\) on \(X\), we can define the corresponding \emph{pushforward measure} \(f_*\nu\) on \(Y\) by letting \((f_*\nu)(A) = \nu(f^{-1}(A))\). It is easy to check that if \(\nu\) is \(G\)-invariant, then \((f_*\nu)\) is also \(G\)-invariant. Moreover,
if \(\nu\) is a probability measure, so is \(f_*\nu\).

If \(G\) is a lcsc group, then there exists a (nonzero) Borel measure \(m\) on \(G\) which
is invariant under the action of \(G\) on itself via left-translations. Such a measure
is called a \emph{Haar measure} on \(G\). If \(m'\) is a second Haar measure on \(G\), then there
exists a positive real constant \(c\) such that \(m'(A)=c m(A)\) for every Borel subset \(A\subseteq G\).
It is well-known that \(G\) is compact if and only if \(m(G)<\infty\). In this case, there exists a unique
Haar probability measure on \(G\). The group \(G\) is said to be \emph{unimodular} if each Haar
measure is also invariant under the action of \(G\) on itself via right translations. It is
well-known that if \(G\) is either discrete, compact, or a connected simple Lie group, then
\(G\) is unimodular.

\begin{proposition}
\label{prop:denseuniquelyergodic}
Suppose that \(\Lambda\) is a countable (discrete) group,  \(G\) is compact, and that \(\Lambda \) is dense in \(G\). Then, the translation action of \(\Lambda \) on \(G\) is uniquely ergodic.
\end{proposition}
\begin{proof}
The probability Haar measure of \(G\) is the unique \(\Lambda\)-invariant probability measure. The details can be found in \cite[Theorem~7.1]{Ada02}.
\end{proof}

Let $E, F$ be countable Borel equivalence relations on $X$ and $Y$ respectively, and let $\mu$ be an $E$-invariant probability measure on $X$. We say that $E$ is \textit{$F$-ergodic} if for every Borel function $f\colon X \to Y$ such that \(x_1\mathbin{E} x_2\implies {f(x_1)}\mathbin{F}{f(x_2)}\) there is a Borel $Z\subseteq X$ with $\mu(Z)=1$ and $f$ maps $Z$ into a single $F$-class.
Note that this naturally generalizes the notion of ergodicity. Propostion~\ref{prop:ergodic} establishes that the action of \(G\) on \((X,\mu)\) is ergodic in the usual measure theoretic sense if and only if the equivalence relation \(E^X_G\) is \(=_Y\)-ergodic, for any standard Borel space \(Y\).

Or proof will make use of the following lemma.

\begin{lemma}
\label{lem:subgroupactionisergodic}
    Let \((X,\mu)\) be a standard Borel probability space and let the action \(\Gamma\curvearrowright (X,\mu)\) be ergodic.
    If \([\Gamma:\Delta]<\infty\), then
    there is a \(\Delta\)-invariant Borel \(Y\subseteq X\)
    with \(\mu(Y)>0\) such that \(\Delta\curvearrowright (Y, \mu)\) is ergodic.
\end{lemma}

\begin{proof}
    Let \([\Gamma:\Delta] = n\), and fix (left-)coset representatives \(\gamma_1,\dotsc,\gamma_n\). 
    Let \(A\subseteq Y\) be any \(\Delta\)-invariant set with \(\mu(A)>0\). Then, note that
    \[
    \Gamma\cdot A = \bigcup_{i\leq n} {\gamma_i \Delta} \cdot A = \bigcup_{i\leq n} {\gamma_i} \cdot A.
    \]
    Therefore, \(\bigcup_{i\leq n} {\gamma_i} \cdot A\) is \(\Gamma\)-invariant. It follows by ergodicity that \(\mu({\gamma_i} \cdot A) = \mu(A)\geq 1/n\).

    Now let \(Y\subseteq X\) be an \(\Delta\)-invariant set of smallest positive measure.
    For any Borel \(A\subseteq Y\) define \(\nu (A) = \mu(A)/\mu(Y)\).
    It follows that if \(A\subseteq Y\) is \(\Delta\)-invariant then \(\mu(A)>0\), then we have  \(\mu(A) = \mu(Y)\), thus \(\nu(A) =1\).
    So the p.m.p action \(\Delta\curvearrowright (Y,\nu)\) is ergodic.
\end{proof}

\subsection{Discrete subgroups of Lie groups}

Now suppose that \(G\) is unimodular and that \(\Gamma\) is a discrete subgroup of \(G\). By the Effros-Mackey cross section theorem (e.g., see~\cite[Theorem~12.17]{Kec}), there exists a Borel transversal \(T\) for \(G/\Gamma\), i.e., a Borel subset \(T\subseteq G\) that meets every coset in a unique element. Such a set \(T\) is also a known as a \emph{fundamental domain} for \(G/\Gamma\).  We fix such a Borel transversal with \(1_G \in T\). Then we can identify \(T\) with \(G/\Gamma\) by identifying each \(t\in T\) with the coset \(t\Gamma\in G/\Gamma\), and the action of \(G\) on \(G/\Gamma\) induces a corresponding Borel action of \(G\) on \(T\) by
\[
g\cdot t = T\cap gt\Gamma.
\]
Let \(\nu\) be the measure induced on \(G/\Gamma\) by \(\mu\restriction T\) via this identification. Then \(\nu\) is \(G\)-invariant. The measure \(\nu\) is called the Haar measure on \(G/\Gamma\) with respect to \(\mu\). It is possible to prove that \(\nu\) is independent of the choice of the transversal \(T\). Finally, we say that \(\Gamma\) is a \emph{lattice} in \(G\) if \(\nu(G/\Gamma)\) is finite. We say that \(\nu(G/\Gamma)\), also denoted by \(\covol_\nu(\Gamma)\) (or simply \(\covol(\Gamma)\)),  is the \emph{covolume} of \(\Gamma\) in \(G\).
The following property of covolume is particularly useful:

\begin{proposition}
    \label{prop:covol}
    If $\Gamma ' \le \Gamma$ are lattices in $G$, then ${\rm{covol}}(\Gamma')=[\Gamma:\Gamma']{\rm{covol}}(\Gamma)$.
\end{proposition}
\begin{proof}
    Let $T$ be a fundamental domain for $G/\Gamma$. Then 
    $$T'=\bigsqcup_{\gamma \in \Gamma / \Gamma'} (\gamma \cdot T)$$
    is a fundamental domain for $G/\Gamma'$, and 
    \[
        \mu(T')= \sum_{\gamma \in \Gamma/\Gamma'} \mu(\gamma \cdot T)
        = [\Gamma:\Gamma']\mu(T).
    \]
    Note that since $\Gamma'$ and $\Gamma$ are lattices, we have 
    \(0< \mu(T),\ \mu(T') < \infty\),
    so $[\Gamma:\Gamma']$ must be finite and 
    \({\rm{covol}}(\Gamma')=[\Gamma:\Gamma']{\rm{covol}}(\Gamma)\).
\end{proof}

\begin{lemma}
\label{lemma : Aut preserves Haar}
    Let \(G\) be a connected semisimple Lie group. Every automorphism of \(G\) preserves the Haar measure. 
\end{lemma}
\begin{proof}
Fix a (left-invariant) Haar measure \(\mu\) on \(G\). For each \(\varphi\in \Aut(G)\) consider the push-forward measure \(\varphi_*(\mu)\) by \(\varphi_*(\mu)(A) = \mu(\varphi^{-1}(A))\). Since the measure \(\varphi_*(\mu)\) is also \(G\)-invariant, there is a unique constant \(\theta(\varphi)>0\) such that \(\varphi_* (\mu) (A) = \theta(\varphi) \mu(A)\) for all Borel \(A\subseteq G\).

\begin{claim}
\label{claim : Aut hom R+}
  The map  \(\theta\colon \Aut(G) \to \mathbb{R}^*_+\) is a group homomorphism.
\end{claim}
\begin{proof}[Proof of Claim~\ref{claim : Aut hom R+}]
    Let \(\varphi,\psi\in \Aut(G)\). Note that \(\theta (\varphi\psi)\) is the unique \(c>0\) such that
    \((\varphi\psi )_*\mu = c\mu\).
    On the other hand, for any Borel \(A\subseteq G\), we have \((\varphi\psi )_*\mu (A)=\theta(\psi) \theta(\varphi)\mu(A)\).

By uniqueness, it follows that \(\theta(\phi\psi) = \theta(\psi) \theta(\varphi)\). Since \(\mathbb{R}^*_+\) is commutative, Claim~\ref{claim : Aut hom R+} is proved.
\end{proof}

Since semisimple connected Lie groups are unimodular \(\mu\) is invariant under inner automorphisms. This means that \(\theta(\Inn(G)) = \{1\}\). So \(\theta\) factors through a Lie group homomorphism \(\Aut(G)/\Inn(G) \to \mathbb{R}_+^*\). Since \(G\) is connected and semisimple, we have
\([\Aut(G): \Inn(G)]<\infty\). (See~\cite[Section~27.4]{Hum}.)  So, since \(\mathbb{R}_+^*\) has
no nontrivial finite subgroups, it follows that \(\theta(\Aut(G)) = \{1\}\), and so \(\mu\) is invariant
under the action of \(\Aut(G)\).
\end{proof}

The following ridigity theorem was proved by Margulis~\cite{Mar77}. (See also Zimmer~\cite[Page~2]{Zim}.)

\begin{theorem}[Mostow--Margulis rigidity theorem]
    \label{thm:mostowmargulis}
    Let $G_0,\ G_1$ be centerless semisimple Lie groups, and suppose $\Gamma_0 \le G_0$ and $\Gamma_1\le G_1$ are irreducible lattices. Suppose $\mathbb{R}\text{-}rank(G_0)\ge 2$. Then any isomorphism $\pi\colon\Gamma_0 \to \Gamma_1$ extends to a continuous isomorphism $G_0 \to G_1$.
\end{theorem}

\subsection{Algebraic groups}

Throughout this paper, let \(\Omega\) be a fixed algebraically closed field of characteristic
\(0\) which contains \(\mathbb{R}\) and all of the \(p\)-adic fields \(\mathbb{Q}_p\). By an \emph{algebraic group}, we
always mean a subgroup \(H\) of some general linear group \(\GL_n(\Omega )\) which is Zariski
closed in \(\GL_n(\Omega )\). If \(H\) can be defined using polynomials with coefficients in the
subfield \(k\) of \(\Omega\), then we say that \(H\) is an algebraic \(k\)-group. For any subring \(R\subseteq \Omega\) ,
we define
\[
\GL_n(R) = \{(a_{ij}) \in \GL_n(\Omega ) \mid a_{ij} \in R\text{ and }\det(a_{ij})^{-1} \in R\};
\]
and if \(H\leq \GL_n(\Omega ) \) is an algebraic group, then we let
\[
H(R) = H \cap \GL_n(R).
\]
In particular, whenever \(R=k\) is a subfield of \(\Omega\) and \(H\) is an algebraic \(k\)-group, then \(H(k)\) is
the group of all matrices in \(H\) with entries in \(k\). If \(H\) is an algebraic \(k\)-group, then \(H\)
is said to be \emph{\(k\)-simple} if every proper normal \(k\)-subgroup is trivial and \emph{almost \(k\)-simple}
if every proper normal \(k\)-subgroup is finite. The \emph{\(k\)-rank} of such a group, denoted by \(\rank_k G\), is the dimension of a maximal \(k\)-split
torus. For instance, we have \(\rank_k \SL_n = n-1\).

Now suppose that \(k\in\GL_n( \Omega)\) is an algebraic \(k\)-group, where \(k\) is either \(\mathbb{R}\) or a 
finite extension of the field \(\mathbb{Q}_p\) of the \(p\)-adic numbers for some prime number \(p\). In this case,
we always regard \(H(k)\leq \GL_n(k)\) as a topological group with respect to the \emph{Hausdorff topology}, which is the topology obtained by restricting the natural topology on \(k^{n^2}\)
to \(H(k)\). With this topology,  \(H(k)\) is a lcsc group. Any topological notion about  \(H(k)\) will refer to this topology. In
particular, if \(k\) is an \(\mathbb{R}\)-group, then \(H(\mathbb{R})^0\) denotes the connected component of the
identity in the Hausdorff topology.

\subsection{\(S\)-arithmetic groups}
Before introducing \(S\)-arithmetic groups and discussing their properties we recall some basic facts and terminology about valuations. 

A \emph{global field} is a finite extension of either the field of rational numbers or a field of rational functions in one variable over a finite field. Throughout this section, let  \(k\) be a global field.  A \emph{valuation} (or \emph{absolute value}) on \(k\) is a function \(v\colon k\to \mathbb{R}^+\cup\{0\}\) such that
\begin{itemizenew}
\item
\(v(x)=0\) if and only if \(x=0\);
\item
\(v(xy)=v(x)v(y)\) for all \(x,y\in k\);
\item
\(v(x+y)\leq v(x)+v(y)\) for all \(x,y\in k\).
\end{itemizenew}
To avoid trivial valuations, we also assume that \(v(x)\neq 1\) for some \(x\in k\), \(x\neq 0\).
Notice that any valuation \(v\) induces a metric \(d_v\) on \(k\) by setting \(d_v(x,y) = v(x-y)\) for all \(x,y\in k\).  We denote by \(k_v\) the completion of \(k\) relative to this metric.  The field \(k_v\) is called a \emph{local field}. Two valuations \(v\) and \(w\) on the same field are said to be \emph{equivalent} precisely when the corresponding metric \(d_v\) and \(d_w\) induce the same topologies. A \emph{place} of \(k\) is a single equivalence class of valuations on \(k\).

A valuation \(v\) is said to be \emph{non-archimedean}  if and only if it satisfies the ultrametric inequality. Otherwise, we say that \(v\) is \emph{archimedean}. 
Note that a global field \(k\) admits an archimedean valuation if and only if \(\operatorname{char} k=0\).
We shall only consider global field \(k\) with \(\operatorname{char} k =0\),  thus the set of archimedean valuations is non-empty.
A valuation \(v\) is archimedean if and only if \(k_v\) is archimedean.  (A local field \(k\) is said to be \emph{archimedean} if and only if it is isomorphic to \(\mathbb{R}\) or \(\mathbb{C}\).)
Moreover,  notice that \(k_v\) is not archimedean if and only if it is totally disconnected. 

The set \(\{x\in k\mid v(x)\leq 1\}\) is called the \emph{ring of integers} of \(k\). This is the unique maximal compact subring of \(k\). When \(v\) is non-archimedean the ring of integers of \(k_v\) is usually denoted by \(\mathcal{O}_v\). More generally, suppose that \(S\) is a finite set of valuations of \(k\). An element \(x\in k\) is said to be \emph{\(S\)-integral} if and only if \(v(x) \leq 1\) for each non-archimedean valuation \(v\notin S\).  The set of \emph{\(S\)-integral elements} is a subring of \(k\) and is denoted by \(k(S)\).

In this paper we mainly deal with the case
 \(k= \mathbb{Q}\). 
To discuss valuations uniformly, it
we let \(\mathbb{Q}_\infty = \mathbb{R}\).
 
 For \(p=\infty\), let \(v_p(x)=|x|\) be the standard absolute value.  If \(p\) is a prime number, then \(v_p\) is the standard \(p\)-adic valuation. That is,  we set 
 \[
 v_p(x)=
 \begin{cases}
     0 &\text{if \(x=0\)},\\
     p^{-\ell} & \text{if \(x=p^\ell \frac{a}{b}\)  and \(a,b\) are coprime with \(p\)}.
 \end{cases}
 \]

 It is well-known theorem of Ostrowski that every valuation on \(\mathbb{Q}\) is equivalent to \(v_p\) for some \(p\in \{p\mid \text{\(p\) is prime}\}\cup\{\infty\}\). (E.g., see~\cite[Theorem~1]{PlaRap}.)  Clearly, \(v_p\) is archimedean if and only if \(p=\infty\).
An element of \(x\in \mathbb{Q}\) is \(S\)-integral if and only if \(x=\frac{m}{n}\) where \(m\in \mathbb{Z}\) and \(n\in \mathbb{N}\)  and the prime factors of \(n\) are in \(S\).
For instance, if \(S=\{p,\infty\}\), then \(\mathbb{Q}(S)= \mathbb{Z}\big[\frac{1}{p}\big]\).

Now, note that is \(G\) is a connected algebraic \(\mathbb{Q}\)-group,  and \(S\subseteq \{p\mid p\text{ prime}\}\cup \{\infty\}\) is a finite subset with \(\infty \in S\),  then \(G(\mathbb{Q}(S))\) embeds into \(\prod_{p\in S} G(\mathbb{Q}_p)\) as a discrete subgroup.  In fact, \(G(\mathbb{Q}(S))\) is a lattice in \(\prod_{p\in S}G(\mathbb{Q}_p)\).  More generally, we consider the following definition.

\begin{definition}
Let \(k\) be a global field and let \(G\) be a \(k\)-group. Also, let \(S\) be a set of inequivalent valuations on \(k\) containing at least the archimedean valuations. 
A subgroup of \(G\) is said to be \emph{\(S\)-arithmetic} if and only if it is commensurable with \(G(k(S))\).
\end{definition}

Before stating the famous \emph{Margulis normal subgroup theorem} for \(S\)-arithmetic groups, we recall the notion of \(S\)-rank.
Suppose that for any \(s\in S\), \(G_s\) is a connected, semisimple \(k_s\)-group. Whenever \(G = \prod_{s\in S} G_s(k_s)\) as above, we
let \(\rank (G) = \sum_{s\in S} \rank_{k_s} G_s\).

\begin{theorem}[Margulis]
\label{thm : s-arith nst}
Suppose that \(G\) is a connected almost \(k\)-simple \(k\)-group.
Let \(\Gamma\) be an \(S\)-arithmetic subgroup of \(G\) and let \(N\) be a normal subgroup of \(\Gamma\). Suppose \(\rank_S G \geq 2\) and that either \(G\) is connected or \(S\) is finite. Then either \(N\) 
lies in the center of \(G\) or \([\Gamma:N]<\infty\).
\end{theorem}

\subsection{Cocycle Superrigidity}
\begin{definition}
    Let $\Gamma$ and $\Lambda$ be countable groups and let $X$ be a standard Borel $\Gamma$-space with an invariant probability measure $\mu$. A Borel function $\alpha\colon \Gamma\times X \to \Lambda$ is called a \textit{cocycle} if for all $g_1, g_2 \in \Gamma$
$$\alpha(g_1g_2, x)=\alpha(g_1, g_2\cdot x)\alpha(g_2,x)$$
for $\mu$-a.e. $x \in X$.
\end{definition}

One particularly relevant way that cocycles can arise is in association to a homomorphism between orbit equivalence relations: let $Y$ be a standard Borel $\Lambda$-space on which $\Lambda$ acts freely and let $f\colon X\to Y$ be a Borel homomorphism between the orbit equivalence relations $E^X_\Gamma$ and $E_\Lambda^Y$. Then we can define an associated Borel cocycle $\alpha\colon \Gamma\times X \to \Lambda$ by setting 
$$\alpha(g,x)=\text{ the unique } h \in \Lambda \text{ such that } h\cdot f(x)=f(g\cdot x).$$

Note that if $\alpha(g,x)=\alpha(g)$ depends only on the $g$ variable, then $\alpha\colon \Gamma\to \Lambda$ is a group homomorphism. 

If the Borel homomorphism $f\colon X \to Y$ that we start with isn't associated to a cocycle that is also a group homomorphism, it may be possible to use an ``adjusting function" to wiggle $f$ and obtain a new homomorphism $f'$ that is associated with a new (but related) Borel cocycle. The relationship between the two cocycles is made precise in the following definition:

\begin{definition}
    Two Borel cocycles $\alpha,\beta\colon \Gamma\times X \to \Lambda$ are \textit{equivalent} if there exists a Borel function $b\colon X \to \Lambda$ such that for all $g \in \Gamma$,
    $$\beta(g,x)=b(g\cdot x)\alpha(g,x)b(x)^{-1}$$
    for $\mu$-a.e. $x \in X$. 
\end{definition}

One can check that if a Borel homomorphism $f\colon X \to Y$ is associated to the cocycle $\alpha\colon \Gamma\times X \to \Lambda$, then $f'\colon X \to Y$ defined by $f'(x)=b(x)\cdot f(x)$ is associated to the equivalent cocycle $\beta\colon \Gamma\times X \to \Lambda$ defined as above. Note that the adjustment only permutes elements within an orbit --- the orbit of $f'(x)$ is the same as the orbit of $f(x)$ for any $x \in X$.  So, \(f\) is a Borel reduction (respectively, a weak Borel reduction) if and only if \(f'\) is a Borel reduction (respectively, a weak Borel reduction).

\begin{definition}
    Let \((X,\mu)\) be a standard Borel probability \(\Gamma\)-space. We say that the action of \(\Gamma\) on \((X,\mu)\) is \emph{cocycle superrigid with countable targets} if for any countable group \(\Lambda\), any Borel cocycle \(\alpha\colon \Gamma\times X \to \Lambda\) is equivalent to some group homomorphism from \(\Gamma\) to \(\Lambda\).
\end{definition}

Let \(G\) be a lcsc group and let \(\pi\colon G \to U(\mathcal{H})\) be a unitary representation of \(G\) on the separable Hilbert space \(\mathcal{H}\). We say that \(\pi\) \emph{almost admits invariant vectors} if for every \(\epsilon > 0\) and every compact subset \(K \subseteq G\),  there exists a unit
vector \(v \in \mathcal{H}\) such that \(\lVert \pi(g)\cdot v - v \rVert < \epsilon\) for all \(g \in K\).
Then we say that \(G\) has \emph{Kazhdan property \((\mathrm{T})\)} if for every unitary representation \(\pi\) of \(G\),  if \(\pi\) almost admits invariant vectors, then \(\pi\) has a non-zero invariant vector.

\begin{proposition}
\label{prop : S-arithmetic (T)}
    Let \(S\) be a finite non-empty set and for each \(s \in S\) let \(k_s\) be a
local field and \(G_s\) a connected semisimple \(k_s\)-group each containing no
connected normal \(k_s\)-subgroup of \(k_s\)-rank \(1\). Then the group
\(G =\prod_{s\in S} G_s (k_s)
\)
has property \((\mathrm{T})\).
\end{proposition}

\begin{example}
We can find a large supply of property \((\mathrm{T})\) groups among \(S\)-arithmetic groups using Proposition~\ref{prop : S-arithmetic (T)}.
First we recall that if \(\Gamma\) is a lattice in the lcsc group \(G\), then \(\Gamma\) has Property (T) if and only if \(G\)
has Property (T). (See \cite[Theorem~1.7.1]{BekdelVal}.)

    \begin{enumerate-(1)}
        \item 
        In the notation of Proposition~\ref{prop : S-arithmetic (T)} let \(k=\mathbb{Q}\), \(S=\{\infty\}\), and \(G_\infty = \SL_n\), for \(n\geq 3\).
        Since \(\rank_\mathbb{R} \SL_n = n-1\), \(\SL_n(\mathbb{R})\) has property (T). Moreover, by the Borel--Harish-Chandra theorem 
        \(\SL_n(\mathbb{Z})\) is a lattice in \(\SL_n(\mathbb{R})\).
        Therefore, the arithmetic group \(\SL_n(\mathbb{Z})\) has property (T) for \(n\geq 3\).     
        \item  Next let \(S=\{\infty, 5\}\), and \(G_\infty=G_5=\SO_5\). The \(S\)-arithmetic group \(\SO_5(\mathbb{Z}[1/p])\) is a lattice in \(\SO_5(\mathbb{R})\times \SO_5(\mathbb{Q}_5)\) with property \((\mathrm{T}) \) because \(\SO_5\) is absolutely almost simple, \(\rank_\mathbb{R} \SO_5 =0\), and \(\rank_\mathbb{R} \SO_5 =2\).
    \end{enumerate-(1)}
\end{example}

The following cocycle superrigidity theorem due to Popa \cite{Pop07} is one of the main tools that will be used in the main argument:

\begin{theorem}[Popa~2007]
    \label{thm:popa}
    Suppose $G$, $\mathbb{G}$ are countable groups containing a countably infinite Kazhdan group $N$ such that $N \trianglelefteq G \le \mathbb{G}$. Then the shift action of \(G\) on \((2)^{\mathbb{G}}\) is cocycle superrigid with countable targets.
\end{theorem}

\begin{remark}
    Theorem~\ref{thm:popa} has since been generalized for malleable actions~\cite{Pop08} and for the free shift action of inner amenable non-amenable groups~\cite{Tuc20}.
\end{remark}

And the following theorem is the second cocycle superrigidity theorem that we will be using.

\begin{theorem}[Furman~2011; Ioana~2016]
    \label{thm:FurIoa}
    Let \(\Gamma\leq G\) be a dense subgroup of a connected, compact group \(G\). Suppose that \(\Gamma\) has Kazhdan property $\mathrm{(T)}$. Consider the left translation action \(\Gamma \curvearrowright(G, m_{G})\), where \(m_{G}\) is the Haar measure of \(G\). Assume that \(\pi_{1}(G)\), the fundamental group of \(G\), is finite. Let \(\Lambda\) be a countable group and \(\alpha \colon \Gamma \times G \to \Lambda \) be a cocycle.
If any group homomorphism \(\pi_{1}(G) \to \Lambda\) is trivial, then \(\alpha\) is cohomologous to a homomorphism \(\delta \colon \Gamma \to \Lambda\).
\end{theorem}

\section{Borel rigidity}
\label{section : main results}

In this section we discuss the proofs of the main results of this paper.

\subsection{Shift-action}
We first address Theorem~\ref{thm:rigidity shift} and focus on the shift action.
\begin{lemma}
    \label{lemma:cocycle rigidity embedding}
    Let \(G\) be a centerless semisimple Lie group with \(\rank_\mathbb{R} (G)\geq 2\), let \(\Gamma\) be an irreducible lattice in $G$, and let \(\Delta\) be any countable group. If \(f\colon (2)^\Gamma \to (2)^\Delta\) is a weak Borel reduction from \(F_\Gamma\) to \(F_\Delta\), then there are:
    \begin{enumerate-(i)}
        \item 
        \label{item:cocycle rigidity embedding(i)}
        a group embedding \(\beta\colon \Gamma \to \Delta\), and
        \item a weak Borel reduction \(f'\colon (2)^\Gamma \to (2)^\Delta\) such that \(f'(g\cdot x) = \beta(g)f'(x)\) almost everywhere.      
    \end{enumerate-(i)}
    Moreover, if \(f\) is a Borel reduction, so is \(f'\).
\end{lemma} 

\begin{proof}
    Let \(f\colon (2)^\Gamma\to (2)^\Delta\) be a countable-to-one Borel map such that \(x_0\mathbin{F}_\Gamma x_1 \implies f(x_0)\mathbin{F}_\Delta f(x_1) \). Then we define a corresponding Borel cocycle $\alpha\colon \Gamma \times X \to \Lambda$ such that $$\alpha(g,x)\cdot f(x) = f(g\cdot x).$$
    Any irreducible lattice in a semisimple Lie group with \(\rank_\mathbb{R}(G)\geq 2\) has property~(T), and therefore the action \(\Gamma\curvearrowright (2)^\Gamma\) is cocycle superrigid by Theorem~\ref{thm:popa}. So, there is some group homomorphism \(\beta\colon \Gamma \to \Delta\) and some Borel function \(b\colon X\to \Delta\) such that
     \(\beta(g) = b(g\cdot x) \alpha(g,x) b(x)^{-1}\) for \(\mu\)-a.e. \(x\in X\).
    Then let \(f'\colon (2)^\Gamma \to (2)^\Delta\) be the corresponding adjustment of \(f\) defined by \(f '(x) = b(x) \cdot f (x)\). Observe that \(f (x)\) is in the same \(\Delta\)-orbit as \(f'(x)\). So \(f'\) is a weak Borel reduction (resp., Borel reduction) if and only if so is \(f\).
    Also note that on a full measure set we have \(f'(g\cdot x)=\beta(g)\cdot f'(x)\), namely \(\beta\) is the cocycle associated with \(f'\).

    Now let \(K= \ker(\beta)\). By the Margulis Normal Subgroup theorem~\cite{Mar}, $K$ is either finite index in \(\Gamma\) or contained in the center of $G$, which is trivial by the assumption on \(G\). Consider the case $[\Gamma:K]<\infty$. We will show that, in fact, this is impossible.
    By Proposition~\ref{lem:subgroupactionisergodic}, there is an $K$-invariant Borel $Y\subseteq X_0$ with $\mu(Y)>0$ such that the induced $K \curvearrowright (Y,\mu_Y)$ is ergodic, where \(\mu_Y\) is defined by $\mu_Y(A)=\frac{\mu(A)}{\mu(Y)}$ for any Borel $A\subseteq Y$. Then, for any $k \in K$, we have 
    \begin{align*}
        f'(k\cdot x)&=\beta(k)\cdot f'(x) = f'(x).
    \end{align*}
    So $f'$ is $K$-invariant. 
    Then by Proposition~\ref{prop:ergodic} there is $Z\subseteq Y$ with $\mu_Y(Z)=1$ such that $f'(Z)=\{x_0\}$ for some $x_0 \in X$. 
    Because $\mu_Y(Z)=1$, we know $\mu(Z)>0,$ so $\mu(\Gamma\cdot Z)>0$. Then since $\Gamma\cdot Z$ is $\Gamma$-invariant, by ergodicity of $\Gamma \curvearrowright (X, \mu)$, we conclude that $\mu(\Gamma\cdot Z)=1$. 

    Now we have \(f'(\Gamma\cdot Z) = \beta(\Gamma)\cdot f'(Z)=  \beta(\Gamma)\cdot x_0 \subseteq \Delta \cdot x_0\).
    So $f'\colon X \to X$ sends a measure one set into a single $F_\Delta^X$-class. This is not possible. In fact, the preimage through \(f'\) of every $F_\Delta^X$-class is countable and therefore has measure zero because \(\mu\) is non-atomic. It follows that \(K\) is trivial, which means that \(\beta\) is a group embedding.
\end{proof}

The following lemma excludes the existence of an embedding in many interesting cases.

\begin{lemma}
    \label{lem:latticeembedding}
    Let $G$ be a centerless semisimple Lie group with $\rank_\mathbb{R}(G)\ge 2$ and let $\Lambda \leq\Gamma$ be irreducible lattices in $G$ with $[\Gamma:\Lambda]<\infty$. If there is an embedding $\Gamma \to \Lambda$, then \(\Gamma = \Lambda\).
\end{lemma}
\begin{proof}
    Let $\phi\colon\Gamma \to \Lambda$ is an embedding and $\Lambda$ is a finite index subgroup, then $\Lambda$ is not a proper subgroup of $\Gamma$ (i.e. $\Lambda=\Gamma$). Let $\Delta=\phi(\Gamma)$, so $\Gamma \cong \Delta$. By Prasad \cite{Pra}, $\Delta$ must be a lattice as well. 
    
    By Theorem \ref{thm:mostowmargulis}, $\phi$ extends to a continuous group isomorphism $G \to G$. In a connected semisimple group, any automorphism preserves the Haar measure. (See Lemma~\ref{lemma : Aut preserves Haar}.) Therefore, \(\Gamma\) and \(\Delta\) have the same covolume. Then
    \[
    {\rm{covol}}(\Gamma)= {\rm{covol}}(\Delta)
        = [\Lambda:\Delta]{\rm{covol}}(\Lambda) 
        =[\Lambda:\Delta][\Gamma:\Lambda]{\rm{covol}}(\Gamma).
    \]
    Thus $[\Gamma:\Lambda]=1$, so $\Lambda=\Gamma$.
\end{proof}

\begin{proof}[Proof of Theorem~\ref{thm:rigidity shift}]
    Let \(\Gamma, \Delta \) be irreducible lattices in the semisimple Lie groups \(G\) and \(H\), respectively.
    Suppose that \(F_\Gamma \sim_B F_\Delta\). Since \(F_\Delta \leq_B F_\Gamma\), there is an embedding \( \rho\colon\Delta \hookrightarrow \Gamma\) by Lemma~\ref{lemma:cocycle rigidity embedding}\ref{item:cocycle rigidity embedding(i)}.
    Let \(\Lambda = \rho(\Delta)\).
    By Prasad~\cite{Pra}, \(\rho(\Gamma)\) is a lattice in \(G\), therefore \([\Gamma :\Lambda]<\infty\). Mostow--Margulis rigidity theorem yields that \(G\cong H\). 
     Since \(\Lambda \cong \Delta\), it is clear that \(F_\Delta \cong_B F_\Lambda\) and therefore \(F_\Gamma \leq_B  F_\Lambda\).
     By Lemma~4.1, there is an embedding \(\rho'\colon \Gamma \to \Lambda\).
    Now it follows that \(\Lambda = \Gamma\) by Lemma~\ref{lem:latticeembedding}. Therefore, \(\Gamma \cong \Delta\) as desired.
\end{proof}

\begin{remark}
\label{remark : Gaussian}
    The proof of Lemma~\ref{lemma:cocycle rigidity embedding} does not use any special properties of the shift action. In particular, we did not use the fact that the shift action \(\Gamma\curvearrowright ((2)^\Gamma,\mu)\) is strongly mixing.  Theorem~\ref{thm:rigidity shift} extends to malleable actions. In particular, it also holds for countable Borel equivalence relations induced by free Gaussian action. For the definition of Gaussian action, we refer the reader to~\cite{PetSin}. A proof of the fact that Gaussian actions satisfying the sufficient condition are malleable, and therefore cocycle superrigid, can be found in~\cite{Fur07}.
\end{remark}

\subsection{The translation action of \(S\)-arithmetic groups}

\label{sec : translation S-arithmetic}

In this subsection we prove Theorem~\ref{thm:rigidity translation}. Throughout this subsection let \(G,H, \Gamma\), and \(\Delta\) as in the statement of Theorem~\ref{thm:rigidity translation}.

The following theorem is an analog of Lemma~\ref{lemma:cocycle rigidity embedding}. The proof closely follows an argument of Ioana~\cite[Theorem~3.2]{Ioa16}.

\begin{theorem}
\label{thm:translation cocycle rigidity embedding}
Let \(m\) be the probability Haar measure on the compact group \(G(k_\infty)\). Let \(\Gamma \curvearrowright (G(k_\infty),m)\) by translation. Let \(\Lambda\) be a finite-index subgroup of \(\Gamma\). We consider the translation actions of \(\Gamma,\Lambda\) on \(G\), respectively. If \(f\colon G(k_\infty) \to G(k_\infty)\) is a weak Borel reduction
from \(E^{G(k_\infty)}_\Gamma\) to \(E^{G(k_\infty)}_\Lambda\), then there exist
\begin{enumerate-(i)}
\item 
a group embedding \(\beta\colon\Gamma \to \Lambda\),
\item 
a \(\Gamma\)-invariant Borel subset \(X_0\subseteq G(k_\infty)\) with \(m(X_0) =1\), and
\item a Borel map \(b\colon X_0 \to \Lambda\)
such that the “adjusted homomorphism” \(f'(x) = b(x)\cdot f(x)\) satisfies
\[
f'(g\cdot x) = \beta(g)\cdot f'(x)
\]
for all \(g\in \Gamma\) and \(x\in X_0\).
\end{enumerate-(i)}
\end{theorem}

\begin{proof}
    Let \(\alpha\colon \Gamma\times G\to \Lambda\) be the cocycle such that \(f(g\cdot x) = \alpha(g,x)\cdot f(x)\). 
    Let \(\widetilde{G}\) be the (algebraic) universal cover of \(G\) and $p\colon \widetilde{G} \to G$ be the finite-to-one covering map and define $\widetilde{\Gamma}=p^{-1}(\Gamma)$. Note that
    \(\widetilde G(k_\infty)\) is compact and simply connected.
    Moreover, \(\Gamma\) has property (T) by Proposition~\ref{prop : S-arithmetic (T)}, so \(\widetilde \Gamma\) has also property (T) by \cite[Remark~3.5.5(ii)]{BekdelVal}. Moreover \(\widetilde \Gamma\) is dense in \(\widetilde G(k_\infty)\) since the covering map is open and onto.

    Now define the lifted cocycle  \(\widetilde{\alpha}\colon \widetilde{\Gamma}\times \tilde G(k_\infty) \to \Lambda\) by $\widetilde{\alpha}(g,x)=\alpha(p(g), p(x))$. By Theorem~\ref{thm:FurIoa}, the cocycle \(\widetilde \alpha\) is cohomologous to a group homomorphism \(\delta\colon \widetilde \Gamma \to \Lambda\). 
    That is, there is a Borel map $\widetilde{b}\colon \widetilde{G}(k_\infty)\to\Lambda$ such that
    $$\delta(g)=\widetilde{b}(g\cdot x)\widetilde{\alpha}(g,x)\widetilde{b}(x)^{-1}$$
    for all $g \in \widetilde{\Gamma}$ and almost every $x \in \widetilde{G}(k_\infty)$. Note that $\widetilde{\Gamma}\curvearrowright (\widetilde{G}(k_\infty), \widetilde{m})$, where $\widetilde{m}$ is the Haar measure on $\widetilde{G}(k_\infty)$, and we can use the adjusting map $\widetilde{b}$ to obtain a Borel homomorphism $\widetilde{f'}\colon \widetilde{G}(k_\infty)\to G(k_\infty)$ from $E_{\widetilde{\Gamma}(k_\infty)}^{\widetilde{G}}$ to $E_\Lambda^{G(k_\infty)}$, defined by $\widetilde{f'}(x)=\widetilde{b}(x)\cdot f(p(x))$. It is easily checked that $\delta$ is the cocycle associated with $\widetilde{f'}$, and we have 
    $$\delta(g)\cdot \widetilde{f'}(x)=\widetilde{f'}(g\cdot x).$$

    Let $r\colon G(k_\infty) \to \widetilde{G}(k_\infty)$ be any Borel retract of $p$. Note that this exists since $p^{-1}$ induces a finite equivalence relation on $\widetilde{G}(k_\infty)$, so there is a Borel selector. Define $\beta\colon \Gamma \to \Lambda$ by $\beta = \delta \circ r$. Note that \(\Lambda\) is \(S\)-arithmetic. --- This follows from Margulis superrgidity theorem for \(S\)-arithmetic groups (e.g., see~\cite[Theorem~10.1.5]{Zim}) and \cite[Theorem~5.9]{PlaRap}.
    
    Before continuing we show that this is well-defined by showing that if $p(g_1)=p(g_2)$ then $\delta(g_1)=\delta(g_2)$. If $p(g_1)=p(g_2)$ then $g_1g_2^{-1} \in \ker(p)$, so $\delta(g_1g_2^{-1}) \in \delta(\ker(p)) \trianglelefteq \delta(\widetilde{\Gamma})$. Using Margulis superrigidity \cite{Zim} and \cite[Theorem~5.9]{PlaRap} we see that $\delta(\widetilde{\Gamma})$ is $S$-arithmetic, so by $\ref{thm : s-arith nst}$ $\delta(\ker(p)) \subseteq Z(G)$, which is trivial. 
    So $\delta(g_1g_2^{-1})=1$. Since $\delta$ is a homomorphism, $\delta(g_1)=\delta(g_2)$. Note that consequently, $\beta$ does not depend on the choice of $r$. 

    We can also show that for almost every $x,x' \in \widetilde{G}(k_\infty)$, if $p(x)=p(x')$ then $\widetilde{b}(x)=\widetilde{b}(x')$. This is because if $p(x)=p(x')$, then we have $x=k\cdot x'$ for some $k \in \ker(p)$. Then since $\widetilde{\alpha}(k,x)=\widetilde{\alpha}(1, x)$ we have $\widetilde{b}(k\cdot x)^{-1}\delta(k)\widetilde{b}(x)=\widetilde{b}(x)^{-1}\delta(1)\widetilde{b}(x)$. Since $\delta$ is a group homomorphism, $\delta(1)=1$, so $\widetilde{b}(k\cdot x)^{-1}\delta(k)\widetilde{b}(x)=1$. But as noted above, $\delta(\ker(p))$ is trivial, so $\delta(k)=1$. Therefore $\widetilde{b}(k\cdot x)=\widetilde{b}(x)$. 

    Now define $b \colon G(k_\infty) \to \Lambda$ by $b = \widetilde{b}\circ r$, and define $f'\colon G(k_\infty) \to G(k_\infty)$ by $f'=\widetilde{f'}\circ r$. We will show that $\beta,\ b$, and $f'$ satisfy the desired properties. 

    First, we see that $\beta$ is injective since $\ker(\beta) \trianglelefteq \Gamma$, so by Theorem~\ref{thm : s-arith nst} $\ker(\beta) \subseteq Z(G)$ which is trivial. 
    
    Since $\widetilde{f'}(x)=\widetilde{b}(x)\cdot f(p(x))$, we have
    \begin{align*}
        f'(x) = \widetilde{f'}(r(x)) = \widetilde{b}(r(x)) \cdot f(p(r(x))) = b(x) \cdot f(x)
    \end{align*}
    as desired.  Additionally, since $f$ is a homomorphism between $E_\Gamma^{G(k_\infty)}$ and $E_\Lambda^{G(k_\infty)}$ and for any $x \in G$, $f'(x) E_\Lambda^{G(k_\infty)} f(x)$ as witnessed by $b(x) \in \Lambda$, $f'$ is also a homomorphism between $E_\Gamma^{G(k_\infty)}$ and $E_\Lambda^{G(k_\infty)}$. Now let $g \in \Gamma$ and $x \in X_0$.
    \begin{align*}
        \beta(g) &= \delta(r(g)) \\
        &= \widetilde{b}(r(g)r(x))\widetilde{\alpha}(r(g),r(x))\widetilde{b}(r(x))^{-1} \\
        \beta(g)\cdot(\widetilde{f'}(r(x))) &= \widetilde{b}(r(g)r(x))\widetilde{\alpha}(r(g),r(x))\widetilde{b}(r(x))^{-1}\widetilde{b}(r(x)) \cdot f(x) \\
        \beta(g)\cdot f'(x) &= \widetilde{b}(r(g)\cdot r(x)) \cdot f(g\cdot x) \\
        b(g\cdot x) \widetilde{b}(r(g)\cdot r(x))^{-1} \beta(g) \cdot f'(x) &= f'(g\cdot x).
    \end{align*}

    Note that $p(r(g\cdot x)) = p((r(g)\cdot r(x)))$, so $\widetilde{b}(r(g\cdot x)) = \widetilde{b}(r(g)\cdot r(x))$. Thus we have $\beta(g)\cdot f'(x) = f'(g\cdot x)$ as desired. 
\end{proof}

The proof of Theorem~\ref{thm:rigidity translation} goes as the proof of Theorem~\ref{thm:rigidity shift}, using the following proposition in place of Lemma~\ref{lem:latticeembedding}.

The following result follows from \cite[Theorem~4.8]{AvnLubMei}.

\begin{proposition}
\label{prop:s-arith cohopfian}
    Let \(\Gamma\) be as in Theorem~\ref{thm:rigidity translation}. Then \(\Gamma\) is co-Hopfian.
\end{proposition}
\begin{proof} The assumption of Theorem~\ref{thm:rigidity translation} implies that \(\rank_S \Gamma\geq 2\). Irreducible higher rank \(S\)-arithemtic groups are cocycle superrigid in the sense of \cite[Definition~4.1]{AvnLubMei}. Then the result follows from \cite[Theorem~4.8]{AvnLubMei}.
\end{proof}

\begin{proof}[Proof of Theorem~\ref{thm:rigidity translation}]
    Let \(G,H,\Gamma, \Delta\) as in the statement of Theorem~\ref{thm:rigidity translation}. If \(E^{G(k_\infty)}_\Gamma \leq_B E^{H(k_\infty)}_\Delta \), then  there is an
embedding \(\beta \colon \Delta \to \Gamma\) as in Theorem~\ref{thm:translation cocycle rigidity embedding}. Let \(\Lambda = \beta(\Delta)\).
Since \(\Lambda \cong \Delta\), \(E^{H(k_\infty)}_\Delta\) is isomorphic to the orbit equivalence relation \(E^{H(k_\infty)}_\Lambda\), induced by a free action \(\Lambda\curvearrowright H(k_\infty)\). If moreover \(E_\Gamma^{G(k_\infty)} \leq_B E^{H(k_\infty)}_\Lambda \), then Theorem~\ref{thm:translation cocycle rigidity embedding} yields an embedding \(\Gamma \to \Lambda\).  It follows from Proposition~\ref{prop:s-arith cohopfian} that \(\Lambda = \Gamma\).
\end{proof}

In addition to the two rigidity theorems discussed in the introduction, we also point out that the same argument as in Theorem~\ref{thm:rigidity shift} and Proposition~\ref{prop:s-arith cohopfian} yield the following:

\begin{theorem}
    Let \(\Gamma\) and \(\Delta\) as in Theorem~\ref{thm:rigidity translation}. Then
    \(F_\Gamma\sim_B F_\Delta\) only if \(\Gamma\cong \Delta\).
\end{theorem}

\section{The containment problem}

In this section we discuss the containment problem; i.e., the problem of finding two countable equivalence relations \(E\subseteq F\) such that \(E\mathbin{\bot}_B F\). This problem was first solved by Adams who provided an example using techniques related to Zimmer's cocycle superrigidity. Thomas~\cite{Tho09} gave an alternative positive answer using the shift action of \(\SL_3(\mathbb{Z})\). We point out that Thomas's argument can be adapted to show the following more general statement:

\begin{theorem}[Thomas~2009]
\label{thm : containment shift}
    Let $G$ be a centerless semisimple Lie group with $\rank_\mathbb{R}(G)\ge 2$ and let $\Lambda \lneq\Gamma$ be lattices in $G$. Then there is $X\subseteq (2)^\Gamma$ with $\mu(X)=1$ such that $E_\Gamma^X \bot_B E_\Lambda^X$.
\end{theorem}

Using material from Section~\ref{sec : translation S-arithmetic} we present an alternative proof of Adams' solution to the containment problem. Our proof is shorter and more conceptual in that it avoids some of the machinery employed by Adams, such  as Ratner's measure classification theorem.

\begin{theorem}
\label{thm : containment translation}
    Let $G$ and $\Gamma$ as in the statement of Theorem~\ref{thm:rigidity translation}, and let $\Lambda$ be a proper finite index subgroup of $\Gamma$. Consider the translation actions of $\Gamma,\ \Lambda$ on $G(k_\infty)$. Then  $E_\Gamma^{G(k_\infty)} \mathbin{\bot}_B E_\Lambda^{G(k_\infty)}$.
\end{theorem}
\begin{proof}
     By Proposition~\ref{prop:s-arith cohopfian}, there is no embedding $\beta\colon \Gamma \to \Lambda$, so $E_\Gamma^{G(k_\infty)} \nleq_B E_\Lambda^{G(k_\infty)}$ because of Theorem~\ref{thm:translation cocycle rigidity embedding}. 
    Towards contradiction, suppose $E_\Lambda^{G(k_\infty)}\le_B E_\Gamma^{G(k_\infty)}$. Let \(m\) be the unique Haar probability measure on \(G(k_\infty)\). By Theorem~\ref{thm:translation cocycle rigidity embedding}, there is
     a \(\Lambda\)-invariant Borel subset \(X \subseteq G(k_\infty)\) with \(m(X) = 1\),
    a group embedding $\beta\colon \Lambda \to\Gamma$ and a Borel reduction $f'\colon X\to X$ from $E_\Lambda^{X}$ to $E_\Gamma^{X}$ such that $f'(h\cdot x)=\beta(h)f'(x)$ for all $h \in \Lambda$ and $x \in X$. 

    We use the following claim:

    \begin{claim}
    \label{claim : full measure image}
        \(m(f'(X)) =1\)
    \end{claim}

    \begin{proof}[Proof of Claim~\ref{claim : full measure image}]
    First, we show that $f'$ is one-to-one. Suppose $f'(x)=f'(y)$. Since \(f'\) is a Borel reduction, we have $x \mathbin{E_\Lambda^X} y$. So there is $h\in \Lambda$ such that $h\cdot x = y$. Then we have
    \[
        \beta(h)\cdot f'(x)=f'(h\cdot x)
        =f'(y)
        = f'(x).\]
    
    Thus $\beta(h)=1$, so by injectivity of $\beta$, $h=1$ and thus $x=y$.

    Now can define a $\beta(\Lambda)$-invariant probability measure $\nu$ on $X$ by pushing forward the probability Haar measure \(m\) on \(G(k_\infty)\). That is, let  $\nu(A)=m(f'^{-1}(A))$ for any Borel \(A\subseteq X\).
    
    Note that the action \(\beta(\Lambda)\curvearrowright (X,\nu)\) is uniquely ergodic. To see this let
    \(H=\overline{\beta(\Lambda)}\) be the closure of the group \(\beta(\Lambda)\) in \(G(k_\infty)\). Then \([G(k_\infty):H]   \leq [\Gamma:\beta(\Lambda)] \). Therefore, \([G(k_\infty):H]\) must be finite.
    On the other hand \(G(k_\infty)\) is connected, so it has no proper subgroup of finite index. It follows that \(H=G(k_\infty)\), thus \(\beta(\Lambda)\) is dense in the compact group \(G(k_\infty)\).
    
    It follows that $\nu=m$, so 
    \[
    m(f'(X))=\nu(f'(X))
        = m(f'^{-1}(f'(X)))
        = m(X)
        = 1.\qedhere
    \]
    \end{proof}

    Continuing the proof of Theorem~\ref{thm : containment translation}, since $m(f'(X))=1$, we have that for all $g\in \Gamma$, $m(f'(X)\cap g\cdot f'(X))=1$. In particular, this intersection is nonempty, so for any fixed $g\in \Gamma$ we can find $x, y \in X$ such that $f'(x)=g\cdot f'(y) \in f'(X)\cap g\cdot f'(X)$. Then $f'(x) E_\Gamma^{X} f'(y)$, so $x E_\Lambda^{X} y$. Thus there is $h \in \Lambda$ where $x=h\cdot y$. Then
    \[
        \beta(h)\cdot f'(y) = f'(h\cdot y)
        = f'(x)
        = g\cdot f'(y).
    \]
    Since $\Gamma$ acts freely on $X$, we have that $g=\beta(h)$. Since $g$ was arbitrary, this shows that $\beta$ is not only an embedding but also an isomorphism. However, $\Gamma \not\cong \Lambda$ by Proposition~\ref{prop:s-arith cohopfian}, so we have a contradiction. Thus $E_\Lambda^{G(k_\infty)}\not\le_B E_\Gamma^{G(k_\infty)}$.
\end{proof}

\end{document}